\documentclass[twocolumn,aps,floatfix,superscriptaddress,longbibliography]{revtex4-2}
\pdfoutput=1 
\usepackage{amsmath,amssymb,eucal,graphicx,float,epstopdf,xparse}
\usepackage{epsfig,subfigure}
\usepackage[utf8]{inputenc}

\usepackage[colorlinks=true, urlcolor=blue, anchorcolor=blue, citecolor=blue,filecolor=blue,linkcolor=blue,menucolor=blue]{hyperref}

\begin{document}

\title{Randomly growing uniform hypergraphs}

\author{P. L. Krapivsky}
\affiliation{Department of Physics, Boston University, Boston, Massachusetts 02215, USA}
\affiliation{Santa Fe Institute, Santa Fe, New Mexico 87501, USA}

\begin{abstract}
We investigate growing uniform random hypergraphs whose hyperedges are $(k-1)$-simplices. At each step, one of two processes occurs: (i) a new vertex is added with probability $p$; (ii) a simplex $\langle v_1,\ldots,v_k\rangle$ with randomly chosen distinct vertices is added with probability $1-p$.  Emerging random uniform hypergraphs exhibit a percolation transition. A giant component is formed when $p<p_c=1-(2k-1)^{-2}$. The fraction of vertices in the giant component scales as $\exp[-C/\sqrt{p_c-p}]$, a signature of the Berezinskii-Kosterlitz-Thouless phase transition. In the non-percolating phase, the cluster size distribution has an algebraic tail with a continuously varying exponent, while at $p=p_c$, it acquires a logarithmic correction. We analyze triadic hypergraphs ($k=3$) in detail and then generalize to higher $k$. We show that the Euler characteristic and the first Betti number of the giant component scale linearly with system size. Higher Betti numbers fluctuate between realizations but remain finite, and hence the normalized Euler characteristic and the first Betti number of the giant component are asymptotically equal, up to a sign. 
\end{abstract}

\maketitle

\section{Introduction} 

Growing random graphs \cite{Hofstad,Newman-book,DM22} constitute an active area of research driven by intellectual challenges and applications \cite{PA-metabolic,PA-international,Barabasi}. Simple stochastic growth rules can produce random graphs having numerous macrohubs, undergoing phase transitions, and exhibiting other striking features \cite{KR14,KR17,KR26}. Some of these models \cite{BKT-Strogatz, KR01-Rodgers,BKT-Sergey,KR02-protein,KR02-BKT,BKT-Bauer-a,BKT-Bauer-b,BKT-Derrida,BK07,Korea18} undergo an infinite-order percolation transition resembling the Berezinskii-Kosterlitz-Thouless (BKT) phase transition found in the two-dimensional spin models \cite{Berezinskii71,Berezinskii72,KT73,Kosterlitz74,Jose,Kosterlitz16} and many other systems \cite{KT69,KT88,KT-jump,KT77,KT-GROSS,KT-Kazakov,KT-Igor}. 

Graphs model systems with binary (pairwise) interactions, while hypergraphs \cite{HG-book} mimic systems with higher-order interactions \cite{Shamir,Newman09,Ginestra18,Ginestra-higher,Bick23}. Edges in graphs, i.e., pairs $\langle v_1 v_2\rangle$ of vertices, are replaced by hyperedges $\langle v_1\cdots v_k\rangle$.  A hyperedge containing $k$ vertices can be represented geometrically by a $(k-1)$-simplex. A hypergraph, however, need not be a simplicial complex because its hyperedge set need not be downward closed.

The simplicial complex induced by a hypergraph is obtained by adding all nonempty subsets of its hyperedges. Simplicial complexes are hypergraphs whose hyperedge sets satisfy the downward-closure property, allowing the application of algebraic topology \cite{Hatcher,Homology}. For a general hypergraph, one can construct an associated simplicial complex and use algebraic-topological tools---such as the Euler characteristic, Betti numbers, and persistent homology---to characterize its structure. These topological characteristics shed light on systems modelled by hypergraphs, e.g., neural, coauthorship, and brain networks, social contagion, etc. \cite{Homology09, Homology19, Brain19, Betti-21b, Petri22a, YZ23}. 

One of the simplest stochastic hypergraph growth processes, the random recursive hypergraph model \cite{KR05, Vazquez22, PK-hyper-23} generates hypergraphs of indefinitely growing dimensionality. Hypergraphs of bounded dimension \cite{Ginestra18, Ginestra-higher} are preferable in many applications.  Here we consider a class of growing random hypergraphs in which all hyperedges have the same cardinality, so the hypergraphs are uniform \cite{HG-book}. We add new vertices and hyperedges independently. Each new hyperedge consists of $k$ distinct vertices selected uniformly at random from the existing vertices. We thus have a class of growing random hypergraphs, GRH($k$), with $k\geq 2$. The GRH(2) process generates graphs, and it was introduced long ago \cite{BKT-Strogatz}. The GRH(3) process is the simplest member of the family that generates genuine hypergraphs. We first study this process in detail (we denote it simply by GRH), and then turn to GRH($k$) processes with arbitrary $k$. 

A uniform triadic hypergraph generated by the GRH process is thus a pair $(\mathcal{V},\mathcal{T})$, where $\mathcal{V}$ is the set of vertices and $\mathcal{T}$ is the set of triangles. By definition, the associated two-dimensional simplicial complex is $(\mathcal{V}, \mathcal{E}, \mathcal{T})$, where $\mathcal{E}$ is the set of all distinct edges of triangles. Ignoring triangles yields the graph $(\mathcal{V}, \mathcal{E})$ associated with the hypergraph.

Two vertices $u$ and $v$ are connected if there is a sequence of triangles $\tau_1, \ldots,\tau_n$ such that $u\in \tau_1, v\in \tau_n$ and adjacent triangles overlap, $\tau_i\cap \tau_{i+1}\ne \emptyset$. This connectivity concept coincides with the graph-theoretic connectivity through a sequence of edges in the graph $(\mathcal{V}, \mathcal{E})$ associated with the hypergraph. By definition, clusters are maximal connected subgraphs of the graph $(\mathcal{V}, \mathcal{E})$. 

The GRH process depends on the parameter $p\in (0,1)$, and it is defined as follows. In each time step, one of the two processes takes place:
\begin{enumerate} 
\item {\bf Vertex creation}.  A new vertex is created with probability $p$.    
\item {\bf Triangle creation}. A new triangle $\langle u, v, w\rangle$ is created with probability $1-p$. The vertices $u, v, w$ are distinct and chosen uniformly and independently from the existing vertices.
\end{enumerate} 

We are mostly interested in the behavior of large random hypergraphs. Constraining the vertices to be distinct does not affect the leading asymptotic behavior of basic observables such as the total number of clusters. Indeed, in the leading order, such observables scale extensively with the number of vertices, and only the sub-leading contribution is affected by the constraint. However, choosing distinct vertices ensures that the induced simplicial complex strictly satisfies the definition of a simplicial complex. We further postulate that once chosen, the same three vertices cannot be chosen again. This requirement is also asymptotically negligible, but it makes the topological analysis cleaner. 

The GRH model exhibits different behaviors depending on whether $p$ is smaller or larger than $p_c=24/25$. When $p_c<p<1$, the cluster size distribution has an algebraic tail with a continuously varying exponent. The size $S$ of the largest cluster diverges algebraically with the number of vertices $M=|\mathcal{V}|$, but remains subextensive
\begin{equation}
\label{S-max}
S \sim M^\frac{1-\sqrt{1-24(1-p)/p}}{2}
\end{equation}
The cluster size distribution acquires a logarithmic correction at $p=p_c=24/25$, and $S\sim (\log M)^{-1}\sqrt{M}$. The giant component containing a finite fraction $g$ of all vertices arises when $p<p_c$. We often simply write `giant' instead of giant component, and we term $g$ the mass of the giant. The mass of the giant exhibits a BKT-type essential singularity: 
\begin{equation}
\label{giant-mass-log}
\log g \simeq - \frac{\pi\sqrt{24}}{5}\,\frac{1}{\sqrt{p_c-p}}
\end{equation}
when $0<p_c-p\ll 1$. This continuous phase transition is sometimes called an infinite-order phase transition: Every derivative of $g$ vanishes at the transition.

The BKT transition in physical systems is often called topological because of the role of topological defects. The BKT-type percolation transitions in networks exhibit the same essential singularity and other mathematical similarities. There is no mapping onto the BKT transition in the XY model or other physical systems. There are no analogs of defects. We examine standard topological characteristics, e.g., Betti numbers. 


All GRH($k$) processes exhibit the BKT-type percolation transition at $p_c=1-(2k-1)^{-2}$. Compared with previous studies of BKT-type percolation transitions in networks, a chief novelty of our work is the analysis of topological characteristics. Algebraic topology is well-suited for studying simplicial complexes, so we extract topological characteristics of a hypergraph by applying algebraic topology to a simplicial complex associated with the hypergraph. With a few exceptions, clusters are topologically trivial (contractible to a point) in the non-percolating phase. The giant is topologically nontrivial; we determine its Euler characteristic and Betti numbers for all $k\geq 2$.

The rest of this paper is organized as follows. In Sec.~\ref{sec:basic}, we derive a few basic properties of the GRH process. We show that $p_c=\frac{24}{25}$ separates the percolating and non-percolating phases. Clusters are finite in the non-percolating phase; in the percolating phase, there is a giant containing a finite fraction $g(p)$ of all vertices. The derivations in Secs.~\ref{sec:basic}--\ref{sec:ns} rely on a differential equation for the generating function encapsulating the cluster size distribution.

In Sec.~\ref{sec:topology}, we probe the topological characteristics of the giant and finite components. We first determine the Euler characteristic, the simplest topological characteristic well-defined for simplicial complexes \cite{Hatcher}. The Euler characteristic of finite clusters is $\chi=1$; more precisely, the fraction of clusters with Euler characteristics $\chi \neq 1$ vanishes in the $M= |\mathcal{V}|\to\infty$ limit. The normalized Euler characteristic and the first Betti number of the giant are equal up to a sign; the second Betti number remains finite. 

In Sec.~\ref{sec:gen}, we analyze GRH($k$) processes with arbitrary $k$. We derive several structural characteristics, such as the cluster size distribution and the mass of the giant, in the same way as for the GRH. The topological characteristics are the least explored even for the GRH(2), the well-studied model \cite{BKT-Strogatz, BKT-Sergey, BK07} with hyperedges being ordinary edges, so the outcome is a growing random graph. We determine the Euler characteristic of the giant component, show that almost all clusters are trees and that the exceptions are unicycles. We analyze the size distribution of unicycles. The tail of this distribution is inversely proportional to the size for all $p\geq p_c$. Therefore, the average number of unicycles scales logarithmically with $M$ when $p\geq p_c$. 
For GRH($k$) models with arbitrary $k$, we show that the normalized Euler characteristic and the first Betti number of the giant are equal up to a sign and given by 
\begin{equation}
\label{Euler-Betti:k}
\frac{\chi_\text{giant}}{M} = \frac{1-p}{p}\,\left[1-kg-(1-g)^k\right] = -\frac{B_1^{(\text{g})}}{M}
\end{equation}
Higher Betti numbers fluctuate from realization to realization, but remain finite in the $M\to\infty$ limit. 

We finish with concluding remarks (Sec.~\ref{sec:remarks}). 

\begin{figure}[t]
\includegraphics[width=0.48\textwidth]{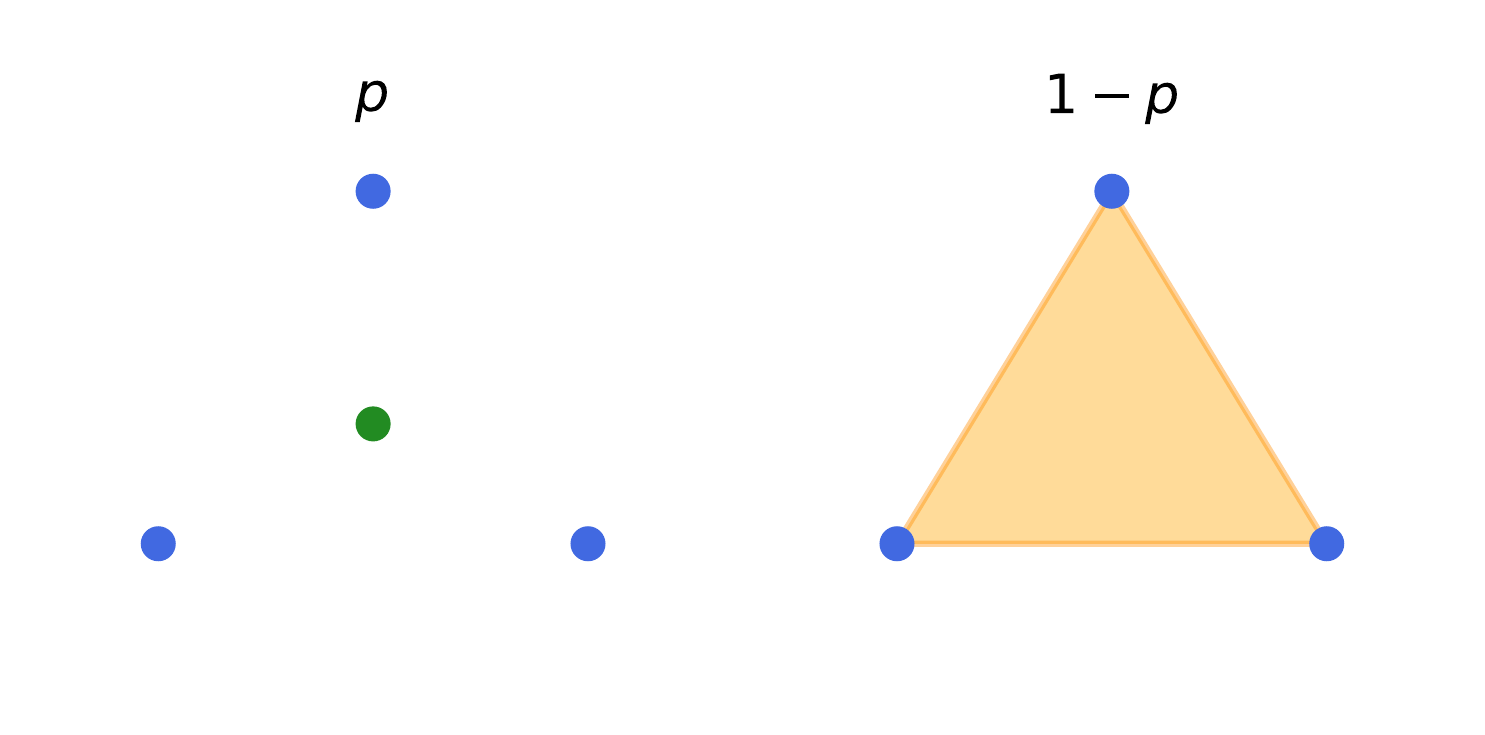}
\caption{An illustration of the GRH process beginning with three vertices (blue circles). With probability $p$, a new vertex (green circle) is added.
With probability $1-p$, a triangle (triadic hyperedge) is added. } 
\label{fig:GRH-process}
\end{figure}

\section{Basic Results}
\label{sec:basic}

In Fig.~\ref{fig:GRH-process}, we visualize the first growth step of the GRH process initialized as a hypergraph with three disjoint vertices. The initial number of vertices should be at least three to allow the choice of three distinct vertices in a triangle creation step. The possibility of a percolation transition is already hinted in the growth law for the total number of clusters following from the rules illustrated in Fig.~\ref{fig:GRH-process}. Let us begin, however, with a simpler quantity $M(t)=|\mathcal{V}(t)|$, the total number of vertices (`$M$' stands for `mass'). The stochastic equation for this quantity 
\begin{equation}
\label{Mt}
M(t+1)=
\begin{cases}
M(t)+1& {\rm prob}~~ p\cr
M(t)    & {\rm prob}~~  1-p
\end{cases} 
\end{equation}
admits an exact solution. In the most interesting large time limit, the total mass is an asymptotically Gaussian random variable with the average and variance 
\begin{equation}
\label{MM}
\langle M\rangle = pt, \quad \langle M^2\rangle - \langle M\rangle^2 = p(1-p)t
\end{equation}

The total number $N(t)$ of clusters seemingly evolves according to the stochastic equation
\begin{equation}
\label{Nt}
N(t+1)=
\begin{cases}
N(t)+1& {\rm prob}~~ p\\
N(t)-2 & {\rm prob}~~1-p
\end{cases} 
\end{equation} 
After adding a vertex, the number of clusters increases by one. This addition process is exactly described by Eq.~\eqref{Nt}. The change in the number of clusters after adding a triangle is more subtle than the change predicted by Eq.~\eqref{Nt}. Namely, the second line in Eq.~\eqref{Nt} becomes erroneous in the percolating phase. The asymptotic growth law, $\langle N\rangle = (3p-2)t$, following from Eq.~\eqref{Nt}, cannot hold when $p\leq \frac{2}{3}$ as the number of clusters is positive: $N(t)\geq 1$. The prediction $(3p-2)t$  is already asymptotically incorrect for $p\leq p_c=\frac{24}{25}$, where a giant is present. 

To correct the second line in Eq.~\eqref{Nt}, we note that the number of clusters remains the same if all three vertices are drawn from the giant, i.e., with probability $g^3$, where $g$ is the fraction of vertices in the giant. The number of clusters decreases by one when two vertices are drawn from the giant; this happens with probability $3g^2(1-g)$.  The number of clusters decreases by two with complementary probability. The corrected version of \eqref{Nt} reads
\begin{equation}
\label{Nt-correct}
N(t+1)=
\begin{cases}
N+1&  {\rm prob}~ p\\
N &   {\rm prob}~(1-p)g^3\\
N-1 &  {\rm prob}~(1-p)3g^2(1-g) \\
N-2 &  {\rm prob}~ (1-p)(1-3g^2+2g^3)
\end{cases} 
\end{equation} 

In the large-time limit, many random variables characterizing the GRH process become sharply concentrated around their mean values. Fluctuations of the total mass are known and given by Eq.~\eqref{MM}. Computing fluctuations in the number of clusters and other basic quantities, e.g., the number of isolated vertices, is much more cumbersome, and little is known even about the variances. However, fluctuations of these random quantities are much smaller than their average values in the $t\to\infty$ limit. Hence, we focus on average values. 

Let $N_s(t)$ be the number of clusters composed of $s$ vertices at time $t$. As stated before, we consider the average numbers but use the same notation since we are interested in the long-time limit where relative fluctuations become asymptotically negligible. In the long-time limit, we can treat time as a continuous variable and replace differences by derivatives: $N_s(t+1)-N_s(t)\to \frac{dN_s}{dt}$. 

The quantities $N_s(t)$ evolve according to 
\begin{eqnarray}
\label{Ns-prod}
\frac{d N_s}{dt} &=& (1-p)\left\{\sum_{i+j+k=s}\frac{ijkN_iN_jN_k}{M^3}
-3\,\frac{sN_s}{M}\right\}\nonumber \\
&+& p\delta_{s,1}
\end{eqnarray}
The terms inside the brackets account for merging of clusters.  The last term on the right-hand side describes the creation of vertices leading to clusters of size one. 

The total number of clusters and the total mass can be expressed via the size distribution:
\begin{subequations}
\begin{align}
\label{N:def}
&N =\sum_{s\geq 1}N_s\\
\label{M:def}
&M =\sum_{s\geq 1} sN_s
\end{align}
\end{subequations}

As a consistency check, we multiply \eqref{Ns-prod} by $s$ and sum over all $s\geq 1$ to deduce $\frac{d M}{dt}=p$, giving $M=pt$, agreeing with the exact average total mass; see \eqref{MM}. Summing \eqref{Ns-prod} over all $s\geq 1$ and using \eqref{M:def} yields $\frac{d N}{dt}=3p-2$ leading to $N=(3p-2)t$. This prediction agrees with \eqref{Nt} and gives the leading asymptotic behavior in the non-percolating $p>p_c$ phase. However, it is clearly incorrect when $p\leq \frac{2}{3}$ as we mentioned before, and even in the entire percolating $p<p_c$ phase. In the percolating phase, the existence of the giant of mass $gM$ implies that the mass of finite components is
\begin{equation}
\label{giant}
\sum_{\text{finite}~s} sN_s(t) = (1-g)M(t)
\end{equation}
and summing Eqs.~\eqref{Ns-prod} over all {\em finite} clusters, we should use \eqref{giant} instead of \eqref{M:def}. Indeed, we obtain
\begin{equation}
\label{N-giant}
\frac{dN}{dt}=p-(1-p)(2-3g^2+g^3)
\end{equation}
The same result follows from \eqref{Nt-correct}. Thus 
\begin{equation}
\label{Nt-av}
\frac{N(t)}{t}=
\begin{cases}
3p-2                            & p_c\leq p\\
3p-2+(1-p)(3g^2-g^3) &p\leq p_c
\end{cases} 
\end{equation} 
when $t\to\infty$. The chief property simplifying the analysis of Eqs.~\eqref{Ns-prod} is extensivity: $N_s(t)$ grow as $M(t)$ for each fixed $s$. Thus we use the ansatz 
\begin{equation}
\label{NsM}
N_s(t)=M(t)\,n_s
\end{equation}
with time-independent $n_s$. 

Inserting \eqref{NsM} into Eq.~\eqref{Ns-prod} yields the recurrence
\begin{equation}
\label{ns-prod}
n_s=\delta_{s,1}+P^{-1}\left\{\frac{1}{3}\sum_{i+j+k=s}ijk\,n_in_j n_k-s\,n_s\,\right\}
\end{equation}
with $P$ defined via 
\begin{equation}
\label{P:def}
P = \frac{p}{3(1-p)}
\end{equation}
Formulas tend to become more compact in terms of $P$, so we mostly use $P$ rather than $p$. Solving Eqs.~\eqref{ns-prod} recurrently one finds that the densities of clusters of even size vanish, $n_2=n_4=\cdots=0$. [The total number of clusters of any even size remains finite: $N_2(t)=0$, $N_4(t)=O(1)$, $N_6(t)=O(1)$, etc.] The non-vanishing densities are
\begin{equation*}
\begin{split}
&n_1 = \frac{P}{1+P}\\
&n_3= \frac{1}{3(3+P)}\left(\frac{P}{1+P}\right)^3\\
&n_5 = \frac{1}{3(5+P)(3+P)}\left(\frac{P}{1+P}\right)^5 
\end{split}
\end{equation*}
etc.  The densities $n_s$ are smooth functions of $P$ over the entire $0<P<\infty$ range. The phase transition at $P=P_c=8$ is not visible from the densities. To detect the phase transition we employ the generating function technique. The generating function 
\begin{equation}
\label{Nz:def}
{\mathcal N}(z)=\sum_{s\geq 1} s\,n_s z^s
\end{equation}
transforms an infinite set of algebraic equations \eqref{ns-prod} into a single nonlinear ordinary differential equation
\begin{equation}
\label{Nz-diff}
z\,\frac{d\mathcal{N}}{dz}= P\,\frac{z-\mathcal{N}}{1-\mathcal{N}^2}
\end{equation}

\begin{figure}[t]
\includegraphics[width=0.48\textwidth]{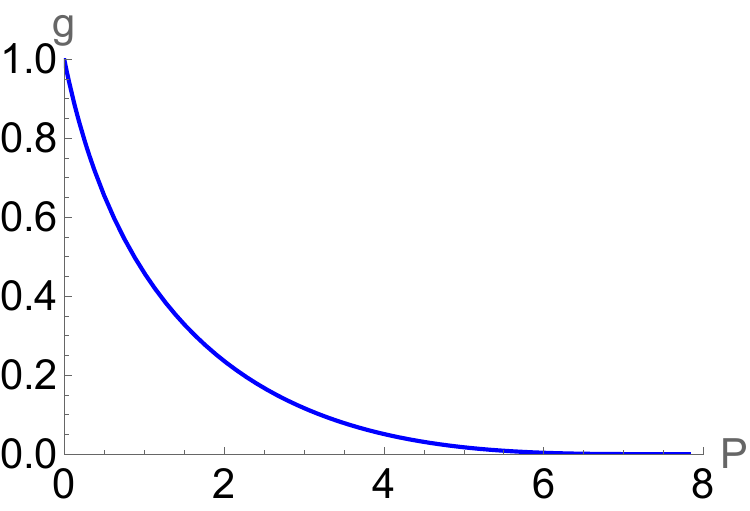}
\caption{The fraction $g$ of vertices in the giant versus $P$ in the percolating regime, $0 < P < 8$. } 
\label{fig:g-GRH}
\end{figure}

To extract the threshold value $P_c=8$ we multiply Eq.~\eqref{Nz-diff} by $(1-\mathcal{N})$, differentiate with respect to $z$, set $z=1$, and recall that in the non-percolating regime
\begin{equation}
\label{N1}
\mathcal{N}(z=1) = \sum_{s\geq 1} s\,n_s = 1
\end{equation}
This procedure leads to the quadratic equation
\begin{align}
\label{m2:P}
2m_2^2 = P(m_2-1)
\end{align}
for the second moment $m_2$ defined via 
\begin{equation}
\label{N2}
\frac{d\mathcal{N}}{dz}\Big|_{z=1} = \sum_{s\geq 1} s^2\,n_s = m_2
\end{equation}

The quadratic equation \eqref{m2:P} has two real roots when $P>8$,  equivalently $p\geq \frac{24}{25}$ when vertex creation is at least 24 times more frequent than triangle creation. The appropriate root 
\begin{equation}
\label{m2}
m_2(P) = \frac{P-\sqrt{P(P-8)}}{4}\,, \qquad P>P_c = 8
\end{equation}
is fixed by the requirement $m_2(\infty)=1$ following from the cluster size distribution $n_s=\delta_{s,1}$ characterizing the model without creation of triangles: $p=1$, or $P=\infty$. 

We first determine the fraction $g$ of vertices in the giant numerically; the result is shown in Fig.~\ref{fig:g-GRH}. The substitution ${\mathcal N}(z)=z F(z)$ reduces \eqref{Nz-diff} to
\begin{equation}
\label{F:eq}
z\,\frac{d F}{dz} + F= P\,\frac{1-F}{1-z^2 F^2}
\end{equation}
Since $F(z)=\sum_{s\geq 1} s\,n_s z^{s-1}$, we have
\begin{equation}
\label{F:IC}
F(0) = n_1 = \frac{P}{1+P}
\end{equation}
We numerically solve Eq.~\eqref{F:eq} subject to the initial condition \eqref{F:IC} and use the identity 
\begin{equation}
\label{N1-super}
F(1) = \mathcal{N}(1) = \sum_{s\geq 1} s\,n_s = 1-g
\end{equation}
following from \eqref{giant} and \eqref{NsM} to determine $g$. 

Using asymptotic techniques, below we extract the behavior near the transition point ($0<P_c-P\ll 1$):
\begin{equation}
\label{giant-mass}
\log g \simeq -\frac{\pi}{\sqrt{1-P/8}}
\end{equation}
The phase transition is of infinite order, or infinitely gentle: Every derivative of $g$ vanishes at the transition.

\begin{figure}[t]
\includegraphics[width=0.48\textwidth]{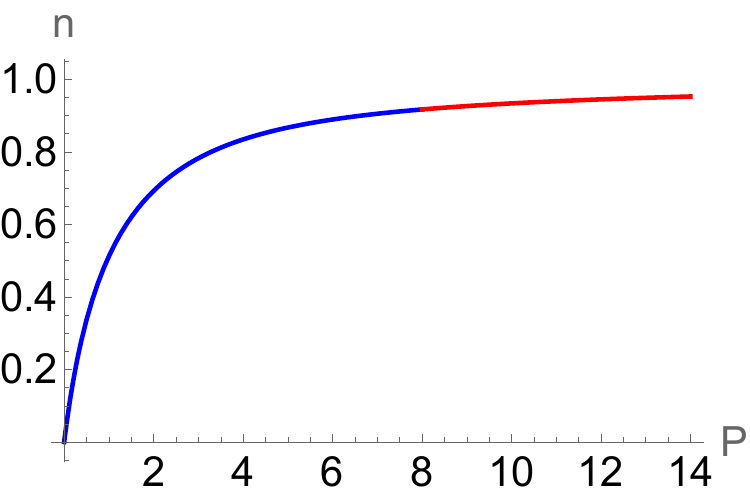}
\caption{The cluster density $n$ given by \eqref{n-prod} as a function of $P$. The cluster density undergoes an infinitely gentle phase transition when $P$ passes through the critical value $P_c=8$.} 
\label{fig:n-GRH}
\end{figure}

Using \eqref{N-giant} we find the cluster density 
\begin{equation}
\label{n-prod}
n=\sum_{s\geq 1}n_s=
\begin{cases}
1-\frac{2}{3P}                                          & 8<P<\infty\\
1-\frac{2}{3P} + \frac{3g^2-g^3}{3P}       & 0<P<8
\end{cases}
\end{equation}
The cluster density depends on the mass $g$ of the giant in the percolating phase; hence, it also exhibits an infinitely gentle phase transition. The transition point is invisible in Fig.~\ref{fig:n-GRH} because all the derivatives on both sides of the transition point are equal.

We already know [cf. \eqref{N1} and \eqref{N1-super}] the mass density, i.e., the first moment of $n_s$:
\begin{equation}
\label{mass-prod}
\sum_{s\geq 1} s n_s=
\begin{cases}
1                                          & 8<P<\infty\\
1-g                                      & 0<P<8
\end{cases}
\end{equation}
The sum in Eq.~\eqref{mass-prod} is over {\em finite} clusters; in the percolating regime, the sum is equal to $1-g$ by the definition of $g$. 
The second moment is 
\begin{equation}
\label{2-prod}
\sum_{s\geq 1} s^2 n_s = \frac{P}{4} \times 
\begin{cases}
1-\sqrt{1-8/P}        & 8<P<\infty\\
\frac{4}{2-g}          & 0<P<8
\end{cases}
\end{equation}
in agreement with \eqref{m2} in the non-percolating regime ($P>8$); the derivation in the percolating regime is presented in Sec.~\ref{subsec:size}. 
The sum in Eq.~\eqref{2-prod} is also over finite clusters; otherwise, the dominant contribution in the $0<P<8$ range would be $g^2 M$. 

\section{Cluster size distribution}
\label{sec:ns}

In classical random graphs, the cluster size distribution is algebraic {\em only} at the critical point; away from criticality, it has exponential tails. In GRH, the cluster size distribution is algebraic in the non-percolating phase; at the critical point, it has an extra logarithmic correction. Therefore, classical random graphs and GRH undergo strikingly different percolation transitions. In this section, we establish these assertions.

\subsection{Non-percolating regime ($P>P_c=8$)}
\label{subsec:sub}

An algebraic $n_s\propto s^{-\gamma}$ tail implies the $(1-z)^{\gamma-2}$ leading singularity \cite{Flajolet} of the generating function \eqref{Nz:def} in the $z\uparrow 1$ limit. To extract the tail we employ the expansion 
\begin{equation}
\label{exp-1}
\mathcal{N}=1+m_2\,(z-1)+\ldots+B(1-z)^{\gamma-2}+\ldots
\end{equation}
The first term in the expansion \eqref{exp-1} follows from normalization \eqref{N1}, while the second regular term in \eqref{exp-1} follows from \eqref{N2}. 
In the expansion  \eqref{exp-1} we also displayed the leading singular term $B(1-z)^{\gamma-2}$. We omitted higher order regular terms $\propto(z-1)^n$ with $n\geq 2$ and higher order singular terms like $(1-z)^{\gamma-1}$. 

\begin{figure}[t]
\includegraphics[width=0.48\textwidth]{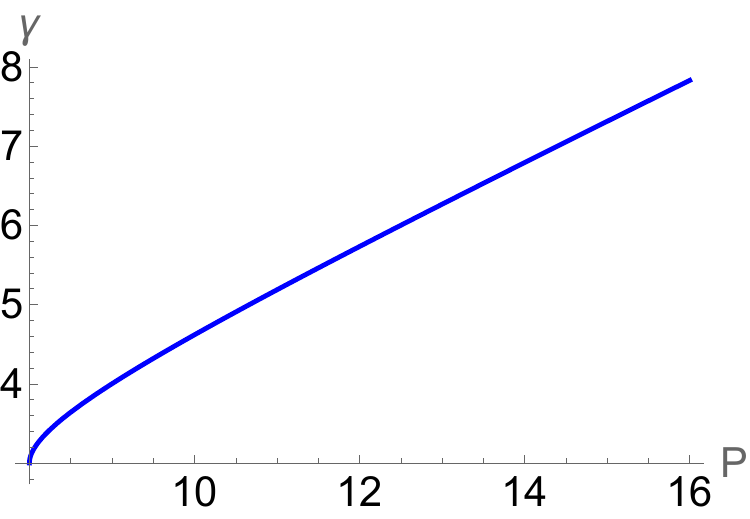}
\caption{The exponent $\gamma$ describing the tail of the cluster size distribution, $n_s\propto s^{-\gamma}$, versus $P$. In the non-percolating and critical regimes, $P\geq P_c=8$, the exponent is given by \eqref{gamma-prod}.} 
\label{fig:gamma3}
\end{figure}

Plugging \eqref{exp-1} into \eqref{Nz-diff} and equating leading regular and singular terms we recover the quadratic equation \eqref{m2:P} for $m_2$ and express the exponent $\gamma$ via $m_2$:
\begin{align}
\label{m2-gamma}
\gamma  = 2 + \frac{1}{m_2-1}
\end{align}
Substituting \eqref{m2} into \eqref{m2-gamma} we fix the exponent 
\begin{equation}
\label{gamma-prod}
n_s\propto s^{-\gamma}, \quad \gamma=2+\frac{4}{P-\sqrt{P(P-8)}-4}
\end{equation}
The asymptotic  behaviors of the exponents are
\begin{equation}
\label{gamma-asymp}
\gamma=
\begin{cases}
3+\sqrt{\frac{P-P_c}{2}} +\frac{P-P_c}{4}+\ldots                                  & P\downarrow 8\\
\frac{P}{2}-\frac{2}{P}+\frac{48}{P^2} -\frac{152}{P^3}+\ldots              & P\uparrow \infty
\end{cases}
\end{equation}
The exponent $\gamma(P)$ increases from $3$ to $\infty$ as $P$ increases from $8$ to $\infty$. Therefore the algebraic decay is steeper than $s^{-3}$ when $P>P_c$, see Fig.~\ref{fig:gamma3}.

\subsection{Critical regime ($P=P_c=8$)}
\label{subsec:crit}

In the critical regime, the exponent predicted by \eqref{gamma-prod} is $\gamma_c=3$. The pure algebraic tail, $n_s\propto s^{-3}$, cannot be correct, however, because it would imply that the sum $\sum_{s\geq 1}s^2n_s$ diverges while Eq.~\eqref{m2} gives the finite value $m_2=\sum_{s\geq 1}s^2n_s=2$ of the sum. The expansion \eqref{exp-1} with $\gamma=3$ has the leading singular term comparable to the second regular term. Therefore, one anticipates the appearance of a logarithmic correction. We thus consider a modified expansion
\begin{equation}
\label{exp-mod}
\mathcal{N}=1+2(z-1)+(z-1)u(z)+\ldots
\end{equation}
where $u(z)$ vanishes slower than any power of $(1-z)$ as $z\uparrow 1$. Plugging \eqref{exp-mod} into \eqref{Nz-diff} yields
\begin{equation*}
(z-1)\,\frac{du}{dz}+\frac{u^2}{2+u}=0
\end{equation*}
which is solved to yield
\begin{equation}
\label{uz}
-\frac{2}{u}+\log |u|+\log(1-z)=\log C
\end{equation}
Finding the integration constant requires integrating the full equation \eqref{Nz-diff}. The leading behavior is independent of the integration constant: $u\to 2/\log(1-z)$ when $z\uparrow 1$. As expected, $u$ vanishes slower than any power of $(1-z)$ as $z\uparrow 1$. Thus
\begin{equation}
\label{exp-crit}
\mathcal{N}=1+2(z-1)+\frac{2(z-1)}{\log(1-z)}+\ldots
\end{equation}
Inverting \eqref{exp-crit} we extract the tail:
\begin{equation}
\label{ns-crit}
n_s\simeq \frac{2}{s^3\,(\log s)^2}\quad {\rm as}\quad s\to\infty.
\end{equation}

\subsection{Percolating regime ($P<P_c=8$)}
\label{subsec:super}

The percolation transition is a continuous transition of an infinite order---all derivatives of $g$ vanish at the transition point. More precisely, near the percolation point, $g$ vanishes according to \eqref{giant-mass}. On the other hand, the behavior of the size distribution is standard---it has an exponential tail. Here, we determine the mass of the giant near the percolation threshold. Mathematically, we need to understand the behavior of the solution of \eqref{Nz-diff} near $z=1$. The asymptotic analysis is possible if, in addition to $1-z\ll 1$, we also assume that $\mathcal{N}$ is close to 1. In the percolating regime, this happens near the percolation transition, see \eqref{N1-super}, where the mass of the giant is small, $g\ll 1$. 

Writing
\begin{equation}
\label{vz-def}
z=1-\zeta, \quad \mathcal{N}(z)=1+\zeta v(\zeta)
\end{equation}
and plugging this into (\ref{Nz-diff}) we obtain in the leading order
\begin{equation*}
2\zeta vv'+2v^2+Pv+P = 0
\end{equation*}
or equivalently
\begin{equation*}
\frac{2v}{2v^2+Pv+P}\,dv+\frac{d\zeta}{\zeta}=0
\end{equation*}
Integrating this equation yields
\begin{subequations}
\begin{eqnarray}
\label{vP:const}
\log(2v^2+Pv+P) + 2 \log \zeta&=& \frac{P}{2a}\,\tan^{-1}\left(\frac{4v+P}{4a}\right)\nonumber \\
&+&{\rm const}
\end{eqnarray}
with 
\begin{equation}
\label{a-def}
a=\tfrac{1}{4}\sqrt{P(8-P)}
\end{equation}
\end{subequations}
From Eqs.~\eqref{N1-super} and \eqref{vz-def} we infer 
\begin{equation}
\label{lim}
\lim_{\zeta\to +0}\zeta v(\zeta)=-g
\end{equation}
Taking the $\zeta\to+0$ limit and using \eqref{lim} we simplify \eqref{vP:const} to give
\begin{eqnarray}
\label{v-zeta}
\log\frac{2v^2+Pv+P}{2(g/\zeta)^2} = \frac{P}{2a}\!\left[\tan^{-1}\left(\frac{4v+P}{4a}\right) + \frac{\pi}{2}\right]
\end{eqnarray}

In the $p\uparrow p_c$ limit we can replace $P$ by $P_c=8$, so that e.g. $2v^2+Pv+P$ is replaced by $2(v+2)^2$. We should keep $a$ and $g$ intact as these quantities go to zero. The leading behavior of the left-hand side of \eqref{v-zeta} is
\begin{equation}
\label{LHS}
2\log(-v-2)+2\log\zeta - 2\log g
\end{equation}
To simplify the right-hand side we use the expansion
\begin{equation}
\label{RHS}
\tan^{-1}\left(\frac{v+2}{a}\right)=\frac{\pi}{2}-\frac{a}{v+2}+\mathcal{O}(a^3)
\end{equation}
Using \eqref{LHS}--\eqref{RHS} we simplify \eqref{v-zeta} to
\begin{equation}
\label{vz-eq}
 -\frac{2}{u} + \log |u| +\log\zeta=\log g +\frac{2\pi}{a}
\end{equation}
where $u=-2-v$. The left-hand sides of Eqs.~\eqref{uz} and \eqref{vz-eq} coincide. Therefore
\begin{equation}
\label{G-eq}
\log g+\frac{2\pi}{a}=\log C
\end{equation}
or equivalently $g=C e^{-2\pi/a}$. The computation of the integration constant $C$ apparently requires the solution of Eq.~\eqref{Nz-diff}. The exact expression \eqref{a-def} is also unnecessary since the analysis is valid near $P_c$ where $a\simeq 2\sqrt{1-P/8}$. Hence
\begin{equation}
\label{giant-inf}
g \simeq C\,e^{-\pi/\sqrt{1-P/8}}
\end{equation}
The amplitude $C$ is not fixed by the asymptotic analysis. 

The fraction of vertices $g$ belonging to the giant is akin to an order parameter characterizing phase transitions. Equation \eqref{giant-inf} shows that the phase transition has an infinite order since every derivative of $g$ vanishes at the transition. The functional form of the order parameter, $g\propto \exp[-C/\sqrt{P_c-P}]$, is common to many growing network models, although the constant $C$ depends on the particular model; see \cite{BKT-Strogatz,KR01-Rodgers,BKT-Sergey,KR02-protein,KR02-BKT,BKT-Bauer-a,BKT-Bauer-b,BKT-Derrida,BK07}. In most models, the pre-exponential factor is constant as in Eq.~\eqref{giant-inf}; a pre-factor algebraically varying with $P_c-P$ has been found in some models \cite{BKT-Derrida}. 

Infinite-order phase transitions in growing network models, growing triadic hypergraphs, and the family GRH($k$) models resemble the BKT phase transition \cite{Berezinskii71, Berezinskii72, KT73, Kosterlitz74} originally found in the two-dimensional $XY$ model in statistical physics. No mapping onto the $XY$ model is known. In the latter context, the correlation length $\xi$ diverges when the temperature approaches the critical temperature from above, $\log \xi\sim 1/\sqrt{T-T_c}$. In the $XY$ model, Kosterlitz \cite{Kosterlitz74} first noticed the square-root divergence. The remarkable universality of this dependence is puzzling; the $(p_c-p)^{-\nu}$ dependence with arbitrary $\nu>0$ seems perfectly acceptable, and an infinite-order phase transition with $\nu=1$ appears even more natural. 



\subsection{The average cluster size and the largest cluster}
\label{subsec:size}

The average mass of {\em finite} clusters is defined via
\begin{equation}
\label{av-size}
\mu=\frac{\sum s^2n_s}{\sum sn_s}=\frac{\mathcal{N}'(1)}{\mathcal{N}(1)}
\end{equation}

\begin{figure}[t]
\includegraphics[width=0.48\textwidth]{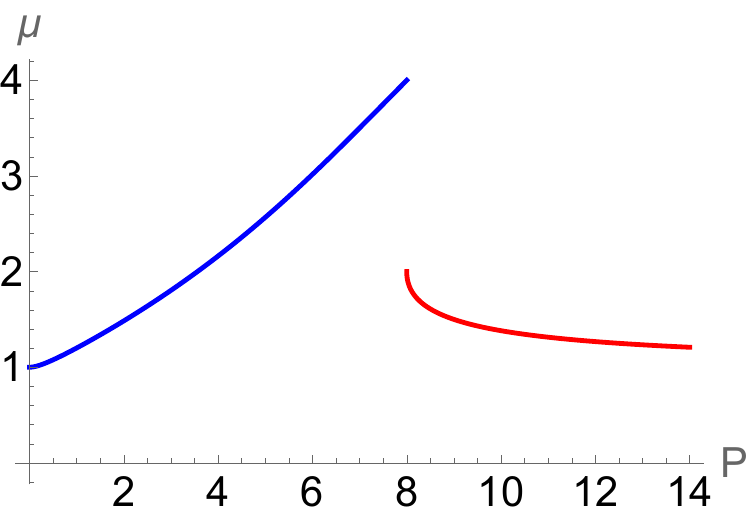}
\caption{The average mass $\mu$ of finite clusters given by \eqref{av-sol} as a function of $P$. The average mass jumps discontinuously at $P=P_c=8$. } 
\label{fig:mu-GRH}
\end{figure}

In the non-percolating regime $\mathcal{N}(1)=1$, $\mathcal{N}'(1)=m_2$, and therefore $\mu=m_2$ with $m_2$ given by \eqref{m2}. In the percolating regime $\mathcal{N}(1)=1-g$, see \eqref{N1-super}. To determine $\mathcal{N}'(1)$ we specialize Eq.~\eqref{Nz-diff} to $z=1$ to give
\begin{equation}
\mathcal{N}'(1) = \frac{P}{1+\mathcal{N}(1)}=\frac{P}{2-g}
\end{equation}
This leads to the announced result \eqref{2-prod}, or equivalently 
\begin{equation}
\label{av-sol}
\mu=\frac{P}{4}\times
\begin{cases}
\frac{4}{(1-g) (2-g)}         &P<8\\
1-\sqrt{1-8/P}                  &P>8
\end{cases}
\end{equation}
This expression is not explicit when $P<P_c$ since we do not have an explicit formula for $g$, but one can plot $\mu(P)$, see Fig.~\ref{fig:mu-GRH}. Since $g\to 0$ as $P\uparrow P_c$,  we know $\mu$ close to the percolation point:
 \begin{equation}
\label{av-jump}
\mu=
\begin{cases}
4       &P\uparrow P_c\\
2       &P\downarrow P_c
\end{cases}
\end{equation}
Although the giant-component fraction undergoes an infinite-order transition, the average mass of finite clusters jumps discontinuously when $P$ crosses the critical point. (Similar jumps occur in other systems exhibiting a BKT phase transition: a jump in the magnetization in Ising chain with long-ranged $r^{-2}$ interactions \cite{KT69,KT88,KT-jump}; a jump in the superfluid density of two-dimensional superfluids \cite{KT77}; etc.)

To estimate the mass $S$ of the largest cluster we employ the criterion $M\sum_{s\geq S}n_s\sim 1$:
\begin{equation}
\label{large}
S\sim 
\begin{cases} 
s^* \log M                           &P<8\\
(\log M)^{-1}\sqrt{M}           &P=8\\
M^\frac{1-\sqrt{1-8/P}}{2}   &P>8
\end{cases}
\end{equation}
The algebraic behavior in the non-percolating phase follows from \eqref{gamma-prod}. The behavior at the critical point is deduced from \eqref{ns-crit}. In the percolating $P<8$ phase, the largest cluster is the giant, of mass $gM$. Equation \eqref{large} instead estimates the mass of the largest finite cluster, which is the second-largest cluster overall. The estimate $S\simeq s^* \log M$ follows from $n_s\propto e^{-s/s^*}$ in the percolating regime. 

\section{Topological characteristics}
\label{sec:topology}

Algebraic topology provides powerful tools for studying topological spaces. The first step is to present a topological space in a manageable form, e.g., as a CW complex or a simplicial complex \cite{Hatcher}. Hypergraphs resemble simplicial complexes. However, one cannot directly apply the methods of algebraic topology to hypergraphs because they are not closed under taking simplicial faces, and downward closure is necessary. One can overcome this drawback by considering a simplicial complex associated with a hypergraph. For a triadic hypergraph, the associated simplicial complex is two-dimensional, $(\mathcal{V}, \mathcal{E}, \mathcal{T})$ for the triadic hypergraph $(\mathcal{V}, \mathcal{T})$. 


In this section, cluster is the shorthand for the two-dimensional simplicial complex representing the cluster.  Each cluster has its own topological characteristics. Below we show that the Euler characteristic and the first Betti number of the giant are extensive. We also compute the Euler characteristic of the entire hypergraph.

\subsection{Euler characteristic}
\label{subsec:Euler-2}

\begin{figure}[t]
\includegraphics[width=0.44\textwidth]{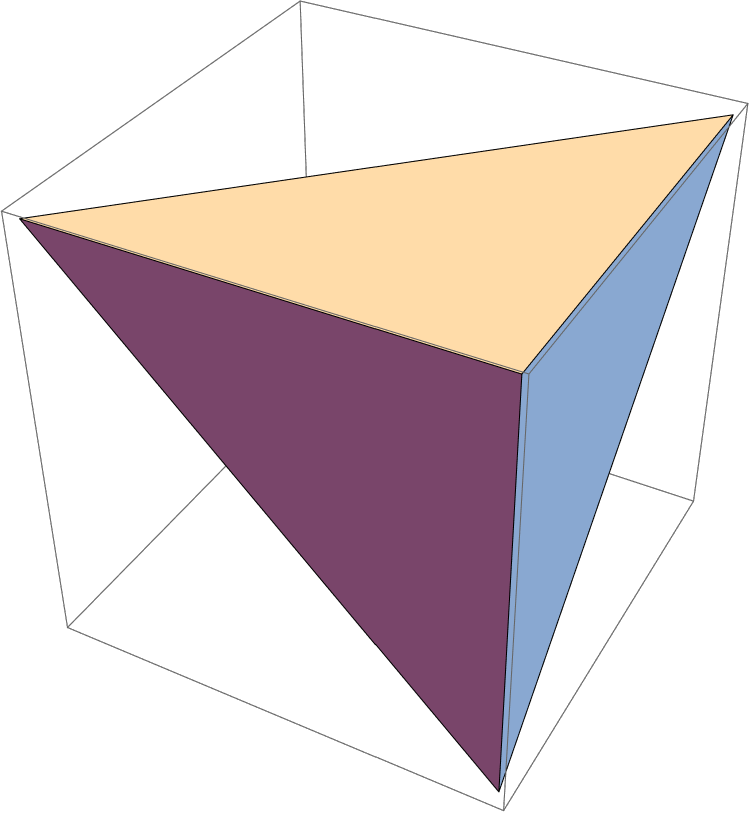}
\caption{The hollow tetrahedron: $(|\mathcal{V}|, |\mathcal{E}|, |\mathcal{T}|)=(4,6,4)$, and hence $\chi=2$. Topologically, the hollow tetrahedron is homeomorphic to the two-dimensional sphere $S^2$. } 
\label{fig:Tetrahedron}
\end{figure}

The Euler characteristic is the simplest topological characteristic of any cluster. We regard clusters as two-dimensional simplicial complexes. Their Euler characteristic is therefore \cite{Hatcher,Dima22}
\begin{equation}
\label{Euler-2}
\chi =  |\mathcal{V}|-|\mathcal{E}|+|\mathcal{T}|
\end{equation}
A finite cluster with $|\mathcal{T}|=s$ triangles has generically $|\mathcal{V}|=2s+1$ vertices and $|\mathcal{E}|=3s$ edges. Thus $\chi=1$ independently of $s$ reflecting that these clusters are topologically trivial, i.e., contractible to a point. This is valid for almost all finite clusters. There are exceptions from the generic decomposition $(\mathcal{V}, \mathcal{E}, \mathcal{T})=(2s+1,3s,s)$. For instance, a cluster with two triangles, $|\mathcal{T}|=2$, and shared edge, the number of vertices is $|\mathcal{V}|=4$ and the number of edges is $|\mathcal{E}|=5$. This cluster is also contractible to a point, and has $\chi=1$. There are also exceptional finite clusters not contractible to a point. The simplest example is the hollow tetrahedron (Fig.~\ref{fig:Tetrahedron}), i.e., a boundary of the tetrahedron which is a three-dimensional simplex. Topologically, the hollow tetrahedron is a sphere, and its Euler characteristic is $\chi=2$. The hollow tetrahedron, as well as finite clusters containing a hollow tetrahedron, can occur with a finite probability during an earlier part of the growth process. Overall, the fraction of exceptional clusters vanishes in the $M\to\infty$ limit.

The topology of the giant is more interesting. The number of vertices evolves according to 
\begin{eqnarray}
\label{V-eq}
\frac{d|\mathcal{V}|}{dt}&=&(1-p)\left[3g^2\sum_{s\geq 1}s^2n_s +3g\sum_{i,j\geq 1}in_i jn_j(i+j)\right] \nonumber \\
&=&(1-p)\left[3g^2 m_2 +6gm_1m_2\right]
\end{eqnarray}
from which
\begin{equation}
\label{V-sol}
\frac{|\mathcal{V}|}{M} = \frac{g^2 m_2 + 2g m_1m_2}{P}
\end{equation}
Recalling that in the percolating phase 
\begin{equation}
\label{m12}
m_1=1-g, \qquad m_2=\frac{P}{2-g}
\end{equation}
we recover the expected result $|\mathcal{V}|=gM$ from \eqref{V-sol}. 

The number of edges in the giant evolves according to 
\begin{equation}
\label{E-eq}
\frac{d|\mathcal{E}|}{dt}=(1-p)\left[3g^3+\frac{9g^2}{2}\,(m_1+ m_2) +9gm_1m_2\right]
\end{equation}
The first term inside the brackets accounts for drawing a triangle within the giant when three edges are added. The second term inside the brackets accounts for drawing a triangle with two vertices inside the giant and one inside a finite component. This term can be written as 
\begin{equation*}
3g^2\sum_{s\geq 1}\left(3\,\frac{s-1}{2}+3\right)sn_s = \frac{9g^2}{2}\,(m_1+ m_2)
\end{equation*}
since a cluster with $s$ vertices is chosen with rate $sn_s$ and such a cluster has $3\frac{s-1}{2}$ edges (since almost all finite clusters are topologically trivial) and three more edges are added. The third term inside the brackets in \eqref{E-eq}  
\begin{equation*}
3g \sum_{i,j\geq 1}in_i jn_j \left(3\,\frac{i-1}{2}+3\,\frac{j-1}{2}+3\right)=9gm_1m_2
\end{equation*}
is derived similarly to the second term. 

\begin{figure}[t]
\includegraphics[width=0.48\textwidth]{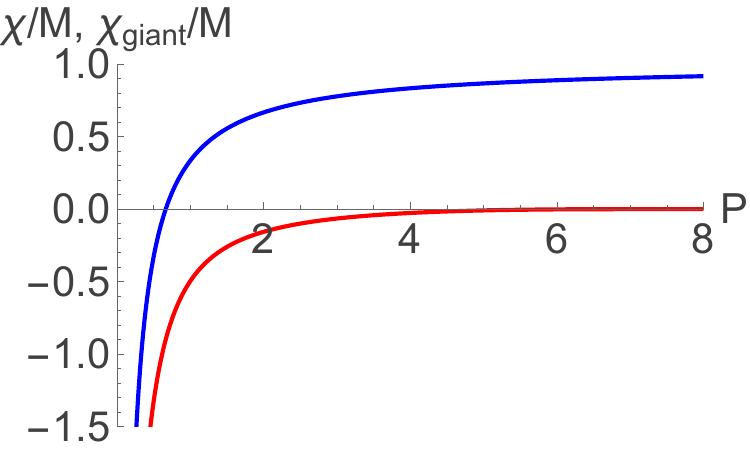}
\caption{Top curve: The normalized Euler characteristic of the entire two-dimensional simplicial complex, Eq.~\eqref{Euler:full}, plotted when $0<P<8$, although it is well-defined for all $0<P<\infty$. Bottom curve: The normalized Euler characteristic of the giant,  Eq.~\eqref{Euler:giant}, plotted in the $0<P<8$ range where the giant exists. The normalized Euler characteristic of the giant is negative, reflecting the extensive number of independent one-dimensional cycles. } 
\label{fig:Euler-char}
\end{figure}

Integrating \eqref{E-eq} and using \eqref{m12} we obtain
\begin{equation}
\label{E-sol}
\frac{|\mathcal{E}|}{M} = \frac{3g}{2}+\frac{3g^2-g^3}{2P}
\end{equation}
The number of edges is three times larger than the number of triangles: $|\mathcal{E}| = 3|\mathcal{T}|$. Hence the normalized Euler characteristic of the giant treated as a two-dimensional simplicial complex reads
\begin{equation}
\label{Euler:giant}
\frac{\chi_\text{giant}}{M} = -\frac{3g^2-g^3}{3P}
\end{equation}

One can also compute the Euler characteristic of the entire hypergraph treated as a two-dimensional simplicial complex. Using $\chi=1$ together with \eqref{n-prod} we find that $n$ is the contribution of the finite components to the normalized Euler characteristic. Thus 
\begin{equation}
\label{Euler:full}
\frac{\chi}{M}= n + \frac{\chi_\text{giant}}{M} = 1-\frac{2}{3P}
\end{equation}
The important point is that the Euler characteristic of the entire complex has a particularly simple expression \eqref{Euler:full} valid for all $0<P<\infty$, even though the topology of the giant is nontrivial.

The normalized Euler characteristic of the giant and the entire hypergraph, Eqs.~\eqref{Euler:giant}--\eqref{Euler:full}, are plotted in Fig.~\ref{fig:Euler-char}. 

\subsection{Betti numbers}
\label{subsec:Betti}

Random simplicial complexes exhibit interesting topological phase transitions \cite{Linial06,Meshulam,Linial16,Kahle,Bobrowski17,Betti-21a,Dima22} often reflected in the behavior of homology groups and the corresponding Betti numbers $B_i$. Intuitively,
\begin{itemize}
\item $B_0$ counts connected components.
\item $B_1$ counts independent one-dimensional cycles (loops).
\item $B_2$ counts independent two-dimensional voids or enclosed surfaces. 
\item  Higher Betti numbers characterize higher-dimensional holes.
\end{itemize}
In the present case of the two-dimensional simplicial complexes, only $B_0, B_1, B_2$ can be nonzero. 

The definition of homology involves an abelian group of coefficients \cite{Hatcher}, integers $\mathbb{Z}$ by default. However, it is simpler to consider 
homology with coefficients in a field $\mathbb{F}$. In this situation, the homology groups of a simplicial complex $\mathcal{K}$ are vector spaces and are therefore characterized, up to isomorphism, by their dimensions $B_k = \text{dim}[H_k(\mathcal{K}; \mathbb{F})]$, known as Betti numbers \cite{Hatcher}, which are independent of the field. 

The Euler characteristic of the giant is equivalently expressed as an alternating sum of the Betti numbers 
\begin{equation}
\label{Euler-2-B}
\chi_\text{giant} =  B_0^{(\text{g})}-B_1^{(\text{g})}+B_2^{(\text{g})}
\end{equation}
The higher Betti numbers vanish. Generally, for the simplicial complex $\mathcal{K}$, the Betti numbers vanish, $B_k(\mathcal{K})=0$ when $k>\text{dim}[\mathcal{K}]$.

The zero Betti number is $B_0^{(\text{g})}=1$ since the giant is connected. The Betti numbers $B_1^{(\text{g})}$ and $B_2^{(\text{g})}$ of the giant are random. Below we argue that the second Betti number varies between realizations but remains finite: it is a non-self-averaging random variable. The first Betti number of the giant is more interesting: It is an asymptotically self-averaging extensive random variable concentrated near its average. Using \eqref{Euler-2-B} together with $B_0^{(\text{g})}=1$ and $B_2^{(\text{g})}=O(1)$, we find 
\begin{equation}
\chi_\text{giant} + B_1^{(\text{g})} =  B_0^{(\text{g})}+B_2^{(\text{g})} = O(1)
\end{equation}
Recalling \eqref{Euler:giant} gives 
\begin{equation}
\label{Betti:giant}
\frac{B_1^{(\text{g})}}{M} = \frac{3g^2-g^3}{3P}
\end{equation}

To prove $B_2^{(\text{g})}=O(1)$ we need some algebraic topology. We give an elementary explanation that requires little knowledge of algebraic topology. We first illustrate how the Betti number $B_2$ can increase during the GRH process. Take a simplicial complex shown in Fig.~\ref{fig:Tetrahedron}, but without the triangle invisible on the figure: $\mathcal{K}=(|\mathcal{V}|, |\mathcal{E}|, |\mathcal{T}|)=(4,6,3)$. This complex is contractible to a point, so its Betti numbers are trivial (apart from $B_0$); in particular, $B_2(\mathcal{K})=0$. If we add the missing triangle to $\mathcal{K}$, turning it into the hollow tetrahedron, we get $B_2=1$. (Indeed, the hollow tetrahedron is homeomorphic to the two-dimensional sphere and $B_2(S^2)=1$.) 

The above observation suggests that $B_2$ can increase by one only by adding a new triangle with {\em all three edges already present}. The underlying idea is that adding a triangle can create a new two-dimensional homology class only when its boundary cycle is already null-homologous in the existing complex. In particular, all three edges must already exist. Adding a triangle with fewer than three edges present does not affect $B_2$. The probability of adding a random triangle with three existing edges vanishes as $t^{-2}$. The sum $\sum_{t\geq 1} t^{-2}$ remains finite, so the number of steps with $B_2(t+1)=B_2(t)+1$ is finite, and $B_2(t)$ converges to a certain finite limit $B_2(\infty)$ which varies from realization to realization. 

Alternatively, one can derive the governing evolution equation for $B_1^{(\text{g})}(t)$ via a kinetic approach repeatedly used before. The first Betti number $B_1^{(\text{g})}(t)$ of the giant can change only when a new triangle is added. One observes:
\begin{enumerate}
\item $B_1^{(\text{g})}(t)$ does not change if three vertices from finite components are used; this happens with probability $(1-g)^3$.
\item $B_1^{(\text{g})}(t)$ does not change if two vertices from finite components and one vertex from the giant are used; this happens with probability $3g(1-g)^2$.
\item $B_1^{(\text{g})}(t)$ increases by one if a vertex from finite components and two vertices from the giant are used; this happens with probability $3g^2(1-g)$.
\item $B_1^{(\text{g})}(t)$ increases by two if three vertices from giant are used; this happens with probability $g^3$.
\end{enumerate}
In the $t\gg 1$ limit, we can ignore exceptional choices, such as choosing two vertices from the same finite component. The above stochastic rules lead to the evolution equation
\begin{equation}
\label{Betti:giant-eq}
\frac{d B_1^{(\text{g})}}{dt} = (1-p)[3g-1+(1-g)^3]
\end{equation}
from which we recover \eqref{Betti:giant}.

\section{Generalized GRH Processes}
\label{sec:gen}

We now consider the family of GRH($k$) processes in which hyperedges are $(k-1)$-simplices with an arbitrary fixed $k\geq 2$. Each new hyperedge is a simplex $\langle u_1,\ldots,u_k\rangle$ with all $k$ vertices chosen uniformly and independently from the existing vertices, with the constraint that they are distinct. As previously, we set the rates of vertex and hyperedge creation steps to $p$ and $1-p$. When $k=3$, we recover the GRH process. The GRH(2) process was introduced in \cite{BKT-Strogatz} and further studied in \cite{BKT-Sergey, BK07}.

Extending many previous results for $k=3$ to the general case is straightforward. In Sec.~\ref{subsec:gen} we present a few general results, mostly those needed for the subsequent analysis of topological characteristics. We then analyze the topological characteristics of the graphs generated by the GRH(2) process in greater depth than the topological analysis of the GRH presented in Sec.~\ref{sec:topology}. For instance, we probe the distribution of unicycles in the non-percolating regime (Sec.~\ref{subsec:uni}). In Sec.~\ref{subsec:Euler}, we return to the GRH($k$) processes with arbitrary $k$, determine the Euler characteristic, and also study the Betti numbers. 

\subsection{General results}
\label{subsec:gen}

Here we present a sample of general results, particularly those that are needed for probing the topological characteristics. 

The recurrence \eqref{ns-prod} generalizes to 
\begin{eqnarray}
\label{ns-k}
Pn_s &=& \frac{1}{k}\sum_{i_1+\ldots+i_k=s} i_1\ldots i_k\,n_{i_1}\ldots n_{i_k}-s\,n_s \nonumber \\
         &+& P \delta_{s,1}
\end{eqnarray}
where we use a more convenient parametrization
\begin{equation}
\label{P-k}
P = \frac{p}{k(1-p)}
\end{equation}
instead of $p$. The generating function \eqref{Nz:def} satisfies an ordinary differential equation
\begin{equation}
\label{Nz-k}
z\,\frac{d\mathcal{N}}{dz}= P\,\frac{z-\mathcal{N}}{1-\mathcal{N}^{k-1}}
\end{equation}

To determine $P_c$ we multiply Eq.~\eqref{Nz-k} by $(1-\mathcal{N})$, differentiate with respect to $z$, set $z=1$, recall that in the non-percolating regime $\mathcal{N}(z=1) = \sum_{s\geq 1} s\,n_s = 1$, and arrive at 
\begin{align}
\label{m2:P-k}
(k-1)m_2^2 = P(m_2-1)
\end{align}
for the second moment $m_2$. The quadratic equation \eqref{m2:P-k} has real roots when $P>P_c=4(k-1)$, equivalently $p>p_c=1-\frac{1}{(2k-1)^2}$. The physically relevant root is
\begin{subequations}
\label{m2-P}
\begin{equation}
\label{m2-k}
m_2=\frac{P-\sqrt{P[P-4(k-1)]}}{2(k-1)}
\end{equation}
In the percolating regime $P<P_c=4(k-1)$, the summation in the second moment $m_2=\mathcal{N}'(1) = \sum_{s\geq 1} s^2 n_s$ is over finite clusters. Specializing \eqref{Nz-k} to $z=1$ and using $\mathcal{N}(1) = 1-g$ we obtain
\begin{equation}
\label{m2-k-perc}
m_2=P\,\frac{g}{1-(1-g)^{k-1}}\,, \quad P < 4(k-1)
\end{equation}
\end{subequations}

Repeating the procedure of Sec.~\ref{subsec:sub} we obtain again \eqref{m2-gamma} fixing the exponent
\begin{equation}
\label{gamma-k}
\gamma = 2 + \frac{2(k-1)}{P-\sqrt{P[P-4(k-1)]}-2(k-1)}
\end{equation}
of the large size tail of the cluster size distribution in the non-percolating regime. 

The critical exponent is universal, $\gamma_c=\gamma(P_c)=3$, for every $k$. Furthermore, in the leading order, the tail is given by \eqref{ns-crit} independently of $k$. Thus, the logarithmic correction and even amplitude are universal. 

In the non-percolating regime, the total number of clusters satisfies the stochastic equation [cf. Eq.~\eqref{Nt}]
\begin{equation}
\label{Nt-k}
N(t+1)=
\begin{cases}
N(t)+1& {\rm prob}~~ p\\
N(t)-k+1 & {\rm prob}~~1-p
\end{cases} 
\end{equation} 
when $t\gg 1$. Therefore the cluster density in the non-percolating regime is 
\begin{equation}
\label{n-k-non}
n = k - \frac{k-1}{p} = 1 - \frac{k-1}{k P} 
\end{equation}

In the percolating regime, $N(t+1)=N(t)+1$ with probability $p$, and  $N(t+1)=N(t)-\ell$ with probability $(1-p)\pi_\ell$, where
\begin{subequations}
\label{pi-ell}
\begin{align}
\pi_\ell = \binom{k}{k-\ell}g^{k-\ell}(1-g)^\ell
\end{align}
with $\ell=0,\ldots,k-2$ is the probability of selecting  $k-\ell$ vertices in the giant. The probability of choosing one vertex in the giant is $\binom{k}{1}g(1-g)^{k-1}=kg(1-g)^{k-1}$, and it leads to $N(t+1)=N(t)-(k-1)$, but the number of clusters also decrease by $k-1$ if we choose all $k$ vertices outside the giant. Thus 
\begin{align}
\pi_{k-1} = k g(1-g)^{k-1} + (1-g)^k
\end{align}
\end{subequations}
Therefore the cluster density in the percolating regime is 
\begin{equation*}
n = 1-\frac{1}{kP}\sum_{\ell=0}^{k-1} \ell \pi_\ell
\end{equation*}
Using \eqref{pi-ell}, we compute the sum and find
\begin{equation}
\label{n-k}
n =  1 - \frac{1-g}{P}  + \frac{(1-g)^k}{k P} 
\end{equation}

\subsection{Topological characteristics of GRH(2)}
\label{subsec:uni}

The GRH(2) process \cite{BKT-Strogatz} is perhaps the simplest random graph process exhibiting the BKT-type percolation transition. However, even for this process, the topological characteristics have not so far been investigated. In the following computations, we use
\begin{subequations}
\begin{align}
\label{n-2}
n &= 
\begin{cases}
1-\frac{1}{2P}                              & P>4\\
 1-\frac{1}{2P} + \frac{g^2}{2P}   & P<4
\end{cases} \\
\label{m2-2}
m_2 &= 
\begin{cases}
\frac{P-\sqrt{P(P-4)}}{2}     & P>4\\
P                                       & P<4
\end{cases}\\
\label{gamma-2}
\gamma &= 1 + \frac{P+\sqrt{P(P-4)}}{2} \qquad (P>4)
\end{align}
\end{subequations}
obtained by specialization of appropriate results from Sec.~\ref{subsec:gen} to $k=2$. 
 
\subsubsection{Euler characteristic of the giant and the entire graph}

The number of vertices in the giant is $|\mathcal{V}|=gM$, while the number of edges increases according to
\begin{equation}
\label{E-1}
\frac{d|\mathcal{E}|}{dt}=(1-p)\left[2g\sum_{s\geq 1}s^2n_s+g^2\right]
\end{equation}
A bond between a giant component and clusters of size $s$ occurs at rate $sn_s$, and since finite components are predominantly trees, we add $s$ edges (a tree of size $s$ has $s-1$ edges). This explains the first term inside the square bracket in \eqref{E-1}. The second term accounts for edge creation inside the giant component. 

Recalling $P=\frac{p}{2(1-p)}$ and $m_2=P$ following from \eqref{m2-k-perc} with $k=2$, we deduce
\begin{equation}
\frac{|\mathcal{E}|}{M} = g+\frac{g^2}{2P}
\end{equation}
from \eqref{E-1}. Hence the normalized Euler characteristic of the giant is 
\begin{equation}
\label{Euler:giant-2}
\frac{\chi_\text{giant}}{M} = \frac{|\mathcal{V}|-|\mathcal{E}|}{M} = -\frac{g^2}{2P}
\end{equation}
Using an alternative formula $B_0-B_1$ for the Euler characteristic, we conclude that $B_1^{(g)}=1-\chi_\text{giant}$. Higher Betti numbers vanish.

Finite components are trees with $\chi(s)=1$; unicyclic components are much more rare as we show in Sec.~\ref{subsub:uni}. Hence, the Euler characteristic of the entire graph is given again by $\frac{\chi}{M}=n + \frac{\chi_\text{giant}}{M}$, cf. Eq.~\eqref{Euler:full}. Using \eqref{Euler:giant-2} together with \eqref{n-2} we find
\begin{equation}
\label{Euler-2-full}
\frac{\chi}{M} = 1-\frac{1}{2P}
\end{equation}
valid for all $0<P<\infty$.

\subsubsection{Unicycles}
\label{subsub:uni}

We now turn to finite components (clusters). We classify the components into three types \cite{Janson93a} depending on the sign of the Euler characteristic: trees are components with Euler characteristic $\chi=1$, unicycles have Euler characteristic $\chi=0$, and components with negative Euler characteristic ($\chi<0$) are called complex. 

In classical random graphs, the total number of unicycles is finite, $\mathfrak{U}=O(1)$, apart from the percolation point where the average number of unicycles diverges logarithmically with the number of vertices 
\begin{equation}
\label{U:CRG}
U \equiv \langle \mathfrak{U}\rangle = \frac{1}{6}\,\log M
\end{equation}
These and more refined results about unicycles are traditionally derived \cite{Stepanov88,Knuth89,Janson93a,Janson93} using probabilistic arguments. One can alternatively derive \eqref{U:CRG} in the realm of the kinetic approach \cite{BK-05}. 

We relied on the kinetic approach in the previous analysis of the cluster size distribution and now adapt it to unicycles in the GRH(2) model. Let $U_s=\langle \mathfrak{U}_s \rangle$  be the average number of unicycles of size $s$. (By size, we mean the number of vertices in a unicyclic component, not the length of a cycle in the component.) The quantities $U_s$ evolve according to 
\begin{equation*}
\frac{d U_s}{dt}=(1-p)\left\{s(s-1)\,\frac{N_s}{M^2}+\sum_{i+j=s}\frac{iU_i jN_j}{M^2}-\frac{sU_s}{M}\right\}
\end{equation*}
which we rewrite as
\begin{equation}
\label{Us-M}
2PM\,\frac{d U_s}{dM} = s(s-1)\,n_s+\sum_{i+j=s}iU_i\,jn_j-sU_s
\end{equation}
The constraint that we choose different vertices has led to the $s(s-1)$ factor.

Let us first analyze the behavior of the average total number of unicycles $U=\sum_{s\geq 2} U_s$. Summing Eqs.~\eqref{Us-M} we find 
\begin{equation}
\label{U-M}
2PM\,\frac{d U}{dM}=m_2-n-Ug
\end{equation}
from which 
\begin{subequations}
\begin{align}
\label{U:GRG}
&U=
\begin{cases}
C_+\log M  & P\geq 4\\
C_-             & P<4
\end{cases} \\
& C_+ = \frac{m_2-n}{2P}\,, \qquad C_- = \frac{m_2-n}{g}
\end{align}
Using Eqs.~\eqref{n-2}--\eqref{m2-2} we obtain 
\begin{align}
\label{Cplus}
C_+ &=  \frac{1-\sqrt{1-4/P}}{4}-\frac{1}{2P}+\frac{1}{4P^2}\\
\label{Cminus}
C_- & = \frac{P-1}{g} + \frac{1-g^2}{2Pg}
\end{align}
\end{subequations}
The growth is logarithmic growth when $P\geq P_c=4$, while in classical random graphs,  the growth is logarithmic only at the percolation point, Eq.~\eqref{U:CRG}. 

We now return to Eqs.~\eqref{Us-M} and note that they can be solved recurrently. The first non-trivial quantity
\begin{equation}
U_2 = n_2 + C M^{-\frac{1}{P}}\,,  \quad  n_2=\frac{1}{2(2+P)}\left(\frac{P}{1+P}\right)^2
\end{equation} 
becomes stationary, $U_2=n_2$, in the $M\to\infty$ limit. The next few exact solutions reveal the same pattern: All $U_s$ approach stationary values in the $M\to\infty$ limit. Thus we seek a stationary (independent of $M$) solution of Eqs.~\eqref{Us-M} and arrive at  
\begin{equation}
\label{Us-rec}
sU_s = s(s-1)\,n_s+\sum_{i+j=s}iU_i\, jn_j
\end{equation}
Introducing the generating function $\mathcal{U}(z) = \sum_{s\geq 1} sU_s z^s$ and using the recurrence  \eqref{Us-rec} we express $\mathcal{U}(z)$ via the generating function \eqref{Nz:def}:
\begin{equation}
\label{U-N}
\mathcal{U} = \frac{z\mathcal{N}' -\mathcal{N}}{1-\mathcal{N}}
\end{equation}
Plugging \eqref{exp-1} into \eqref{U-N}  when $P>P_c$ and \eqref{exp-crit} into \eqref{U-N} when $P=P_c$ we find $\mathcal{U} = \frac{1-\sqrt{1-4/P}}{2}\,\frac{1}{1-z}$ when $z\uparrow 1$ from which we extract the large $s$ behavior 
\begin{equation}
\label{Us-1}
U_s\simeq \frac{1-\sqrt{1-4/P}}{2}\,\,s^{-1}
\end{equation}
for all $P\geq P_c=4$. Therefore, the size distribution of unicycles exhibits remarkably universal behavior, manifested by a universal $s^{-1}$ tail at and above the percolation transition. An algebraic tail of the cluster-size distribution in non-percolating regimes is a distinctive feature of the BKT-type transition. However, the decay exponent $\gamma=\gamma(P)$ depends on $P$. Furthermore, an additional logarithmic correction appears at the transition point; see \eqref{ns-crit}. In contrast, the size distribution of unicycles shows the same universal behavior for all $P\geq P_c$. 

When $0<P_c-P\ll 1$, the unicycle size distribution has an exponential tail
\begin{equation}
\label{Us-star}
U_s\propto e^{-s/s^*}, \quad 
\log s^* \simeq \frac{\pi}{\sqrt{1-P/P_c}}
\end{equation}
The derivation of Eq.~\eqref{Us-star} relies on the properties of the generating function $\mathcal{N}(z)$ in the percolating regime close to the percolation point, $0<P_c-P\ll 1$. The analysis is the same as for the GRH model, see Sec.~\ref{subsec:super}. 

Trees are much more abundant than unicycles. The number of trees of size $s$ is  an asymptotically self-averaging random quantity concentrated around the average $Mn_s$ with fluctuations of order $\sqrt{M}$. The average number $U_s=\langle \mathfrak{U}_s \rangle$ of unicycles of size $s$ is finite. Therefore $\mathfrak{U}_s$ is a non-self-averaging random quantity and the average $U_s$ sheds limited light on the distribution $P_s(u)=\text{Prob}[\mathfrak{U}_s =u]$. 

The total number of unicycles, Eqs.~\eqref{U:GRG}, appears to be an asymptotically self-averaging random quantity. The average grows logarithmically in system size, and fluctuations are expected to scale as $\sqrt{\log M}$. In Appendix~\ref{ap:uni-mass} we analyze the total mass of unicyclic components: $W=\sum_{s\geq 1} sU_s$. This random quantity appears asymptotically self-averaging, growing on average as
\begin{equation}
\label{mass:uni}
W \sim 
\begin{cases}
M^\frac{1-\sqrt{1-4/P}}{4}    & P >\frac{25}{6}  \\
M^\frac{1}{5}\log M              & P = \frac{25}{6} \\
M^\frac{1-3\sqrt{1-4/P}}{2}  & 4<P<\frac{25}{6} \\
\text{finite}                           & P<4
\end{cases}
\end{equation}

It is instructive to compare the average total mass of unicyclic components with the mass $S$ of the largest cluster. Similarly to the derivation of Eq.~\eqref{large} in the case of triadic hypergraphs, we find 
\begin{equation}
\label{large-2}
S\sim 
\begin{cases} 
s^* \log M                           &P<4\\
(\log M)^{-1}\sqrt{M}           &P=4\\
M^\frac{1-\sqrt{1-4/P}}{2}   &P>4
\end{cases}
\end{equation}
The growth exponent in \eqref{large-2} is larger than in \eqref{mass:uni}, exactly twice larger when $P>\frac{25}{6}$.

\subsection{Topological characteristics of GRH($k$)}
\label{subsec:Euler}

We computed the Euler characteristic when $k=3$ and $k=2$. The corresponding expressions, Eqs.~\eqref{Euler:full} and \eqref{Euler-2-full}, suggest the general formula 
\begin{equation}
\label{Euler-k}
\frac{\chi}{M} = 1-\frac{k-1}{kP}
\end{equation}
which we now derive using the general relation 
\begin{equation}
\label{Euler-n}
\frac{\chi}{M}=n + \frac{\chi_\text{giant}}{M}
\end{equation}
Equation \eqref{Euler-n} reflects that almost all finite components are topologically trivial (contractible to a point). Hence, they have the same Euler characteristic, $\chi=1$, independently of the size. Therefore, in the leading order, the contribution of finite components to the normalized Euler characteristic is $\sum_s n_s \chi=n$. We know $n$, Eqs.~\eqref{n-k-non} and \eqref{n-k}, for arbitrary $k$. We now need to compute the normalized Euler characteristic of the giant for arbitrary $k$. 

Consider a hyperedge-creation event and suppose that exactly $k-\ell$ of its $k$ vertices belong to the giant immediately before the event. The remaining $\ell$ vertices belong, with probability tending to one, to distinct finite components. Those components become incorporated into the giant. The Euler characteristic changes by $\Delta \chi_\text{giant}=1-(k-\ell)$. The newly added simplex has Euler characteristic one, and its intersection with the previously existing components consists of  $k-\ell$ isolated vertices in the giant and one vertex in each absorbed finite component. We already know the probability $\pi_\ell$; see \eqref{pi-ell}. Therefore
\begin{eqnarray}
\label{Euler-k:eq}
\frac{d\chi_\text{giant}}{dt} &=& (1-p)\sum_{\ell=0}^{k-1} \binom{k}{k-\ell}g^{k-\ell}(1-g)^\ell (1-k+\ell) \nonumber \\
&=& (1-p)\big[1-kg-(1-g)^k\big]
\end{eqnarray}
from which
\begin{equation}
\label{Euler-giant-k}
\frac{\chi_\text{giant}}{M}=\frac{1-kg-(1-g)^k}{kP}
\end{equation}
in the percolating regime. Combining \eqref{Euler-n}, \eqref{n-k-non}, \eqref{n-k}, and \eqref{Euler-giant-k} we arrive at \eqref{Euler-k}. 

The kinetic equation for first Betti number of the giant is derived similarly to \eqref{Betti:giant-eq} for $k=3$, and to \eqref{Euler-k:eq}:
\begin{eqnarray}
\label{Betti-1-k:eq}
\frac{dB_1^{(\text{g})}}{dt} &=& (1-p)\sum_{\ell=0}^{k-1} \binom{k}{k-\ell}g^{k-\ell}(1-g)^\ell (k-\ell-1) \nonumber \\
&=& (1-p)\big[kg-1+(1-g)^k\big]
\end{eqnarray}
from which 
\begin{equation}
\label{Betti-1-giant-k}
\frac{B_1^{(\text{g})}}{M}=-\frac{\chi_\text{giant}}{M}=\frac{kg-1+(1-g)^k}{kP}
\end{equation}

Higher Betti numbers can grow when a new simplex $\sigma$ is attached to the existing complex $\mathcal{K}$. This can be seen by applying a Mayer-Vietoris sequence to the complex $\mathcal{L} = \sigma \cap \mathcal{K}$. The Betti numbers $B_i(\sigma \cup \mathcal{K})$ with $i \geq 2$ can grow only if the reduced homology group $\widetilde{H}_{i-1}(\mathcal{L})$ is nonzero. This is the key topological observation. (See \cite{Hatcher, Homology} for reduced homology and a Mayer-Vietoris sequence.) The group $\widetilde{H}_{i-1}(\mathcal{L})$ is nontrivial if there is a cycle in the 1-skeleton of $\mathcal{L}$. The probability of such an event can then be bounded by considering the possible cycles among the fixed number $k$ of selected vertices. As in Sec.~\ref{subsec:Betti}, one concludes that this probability decays like $t^{-2}$. Only finitely many such events occur, and hence $B_i = O(1)$ for $i=2,\ldots,k-1$. Using 
\begin{equation}
\label{Euler-k-B}
\chi_\text{giant} =  \sum_{i=0}^{k-1} (-1)^i B_i^{(\text{g})}
\end{equation}
we thus conclude that 
\begin{equation}
\label{Euler-Betti}
\chi_\text{giant} +  B_1^{(\text{g})} =  1+\sum_{i=2}^{k-1} (-1)^i B_i^{(\text{g})} = O(1)
\end{equation}

Near the percolation point
\begin{equation}
\label{B1-giant-k}
\frac{B_1^{(\text{g})}}{M}=-\frac{\chi_\text{giant}}{M} \simeq \frac{k-1}{2P_c}\,g^2
\end{equation}
The essential singularity is therefore twice as strong, in the exponent, as that of the giant-component mass.

\section{Concluding Remarks}
\label{sec:remarks}

A class GRH($k$) of growing random hypergraphs, with hyperedges being $(k-1)$-simplices, exhibits the BKT-type percolation transition for any fixed $k\geq 2$. Analyzing the phase transition is necessary even when the goal is to determine topological characteristics. Here are the main results that we established:
\begin{itemize}
\item The Euler characteristic of the giant is extensive.
\item The first Betti number of the giant is extensive and asymptotically determined by the Euler characteristic.
\item Higher Betti numbers of the giant remain finite, although they fluctuate between realizations. 
\end{itemize}

Betti numbers generally provide much finer topological information than the Euler characteristic. For the giant component of the GRH($k$) processes, however, the two descriptions are asymptotically equivalent: the higher Betti numbers remain finite, while the first Betti number is determined, to leading order, by the Euler characteristic.

The Euler characteristic of the entire hypergraph is also extensive. The normalized Euler characteristic of  the entire hypergraph, Eq.~\eqref{Euler-k}, remains a smooth function over the entire range of the parameter $p$.  

We analyzed not only the topological characteristics of the giant component in GRH(2), but also finite components in the non-percolating regime. Almost all finite components are trees; unicycles (Euler characteristic $\chi=0$) are rare, and complex components with $\chi<0$ are exceedingly rare. Even for GRH(2), we do not know answers to several basic questions. For instance, in evolving classical random graphs, there is never more than a single complex component with probability $5\pi/18$. Up to the percolation point, the classical evolving random graph consists entirely of trees, unicycles, and bicycles with probability $\sqrt{2/3}\cosh\big(\sqrt{5/18}\big)$; see \cite{Stepanov88, Knuth89, Janson93a, Janson93} for these and other subtle behaviors. Developing a comparable understanding for GRH(2) is an intriguing challenge. 

Studying the limiting distributions of higher Betti numbers of the giant in the GRH($k$) processes is left for future work. For the triadic model, there is only one such limiting distribution
\begin{equation}
P(m) = \text{Prob}[B_2^{(\text{g})}=m], \qquad m=0,1,2,\ldots
\end{equation}
The initial condition affects the limiting distribution when $m$ is small; the large $m$ tail is perhaps universal. 

Homology groups are less intuitive than homotopy groups  \cite{Hatcher}, but they are considerably easier to compute \cite{Homology}. It would also be interesting to investigate the homotopy groups, at least the fundamental group of GRH($k$). Among other useful topological invariants, we mention Hodge Laplacians.

A class of models with preferential attachment of vertices during the triangle-creation step was proposed in \cite{Betti-21a}: a triangle $\langle u, v, w\rangle$ is added with probability proportional to $(d_u+a)(d_v+a)(d_w+a)$, where $d_v$ is the facet degree of vertex $v$. Simulations \cite{Betti-21a} indicate that the normalized first Betti number of the giant is again positive in the percolating phase, $p<p_c(a)$; the normalized second Betti number of the giant also becomes positive later, $p<p_*(a)$, with $p_*<p_c$. Simulations also suggest \cite{Betti-21a} that the homological percolation transition resembles the BKT-type transition with a caveat: When $0<p_*(a)-p\ll 1$
\begin{align}
\label{b2}
\log \frac{B_2^{(\text{g})}}{M}  \sim -[p_*(a)-p]^{-\beta(a)} 
\end{align}
with $\beta(a)>\frac{1}{2}$; almost all BKT-type phase transitions are characterized by the universal exponent $\beta=\frac{1}{2}$, cf. Eq.~\eqref{giant-mass}. It would be interesting to prove analytically the appearance of the homological percolation point at $p_*(a)$ which is strictly positive and precedes the percolation transition: $0<p_*(a)<p_c(a)$.

\bigskip\noindent
I am grateful to H. Hartle for numerous useful comments.

\bibliography{references-graphs}

\appendix
\section{Mass of unicycles in the non-percolating phase}
\label{ap:uni-mass}

We define the total average mass, or weight, of unicycles via $W=\sum_{s\geq 1} sU_s$ . In the non-percolating phase, $W$ evolves according to equation 
\begin{equation}
\label{U-mass}
2PM\,\frac{d W}{dM} = m_3-m_2 + m_2 W
\end{equation}
derived by multiplying Eq.~\eqref{Us-M} by $s$, summing over all $s$. (In the percolating phase, $\sum_{s\geq 1} s^2U_s$ appears in the right-hand side.)

The second moment $m_2$ is given by Eq.~\eqref{m2-2} in the non-percolating phase. The third moment of the cluster size distribution, $m_3=\sum_{s\geq 1} s^3 n_s$, remains finite in the thermodynamic limit when $\gamma>4$, that is, $P>9/2$. To determine $m_3$, we use the governing equation for the generating function 
\begin{equation}
\label{N-eq-2}
z(1-\mathcal{N})\frac{d\mathcal{N}}{dz}= P(z-\mathcal{N})
\end{equation}
Differentiating \eqref{N-eq-2} twice with respect to $z$ and setting $z=1$ we arrive at
\begin{align}
\label{m3-m2}
m_3 = \frac{m_2}{3-2m_2}
\end{align}
In the derivation, we have used $m_2^2 = P(m_2-1)$ and $\mathcal{N}''(1)=m_3-m_2$. 

\begin{figure}[t]
\includegraphics[width=0.44\textwidth]{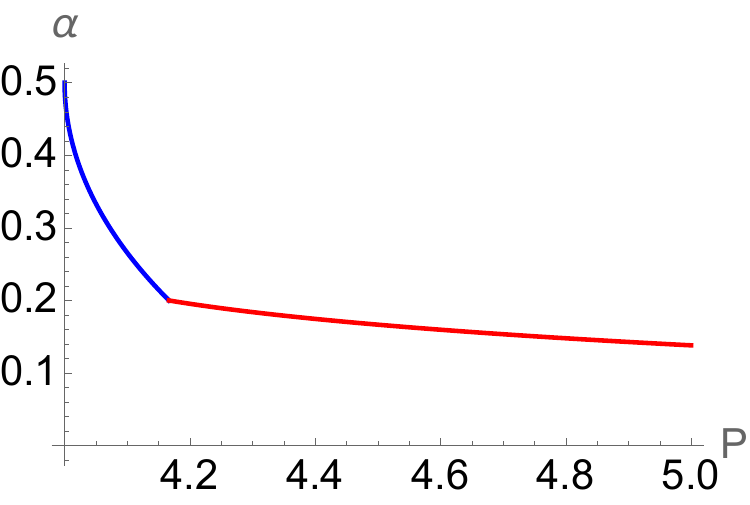}
\caption{The exponent $\alpha$, Eq.~\eqref{alpha:P}, versus parameter $P$. The exponent accounts for the scaling of the mass of unicycles, Eq.~\eqref{alpha:def}. The exponent undergoes a continuous phase transition when $P$ passes through $P=\frac{25}{6}$. The mathematical reason is that the last term on the right-hand side of \eqref{U-mass} dominates in the $P>\frac{25}{6}$ range, while the first term dominates when $4<P<\frac{25}{6}$. } 
\label{fig:alpha}
\end{figure}

Therefore when $P>9/2$, the last term on the right-hand side of \eqref{U-mass} dominates and we find that the average mass of unicycles scales algebraically
\begin{equation}
\label{alpha:def}
W\sim M^\alpha
\end{equation}
with $\alpha=\frac{m_2}{2P}$. If $P<9/2$, the third moment diverges 
\begin{equation}
\label{m3:M}
m_3\sim \sum_{s<S} s^3 n_s\sim S^{4-\gamma}\sim M^\frac{4-\gamma}{\gamma-1}
\end{equation}
since $\gamma<4$. The scaling $S\sim M^\frac{1}{\gamma-1}$ of largest cluster size in the non-percolating regime was found from the extreme statistics criterion: $\sum_{s>S}s^{-\gamma} \sim M^{-1}$. The decay exponent $\gamma$ is given by Eq.~\eqref{gamma-2}. 

Using \eqref{U-mass} and \eqref{m3:M} and analyzing various possibilities we find that the last term on the right-hand side dominates when $P>\frac{25}{6}$ giving again $\alpha=\frac{m_2}{2P}$, while the first dominates when $4<P<\frac{25}{6}$. Thus the average mass of unicycles scales according to \eqref{alpha:def} with exponent 
\begin{equation}
\label{alpha:P}
\alpha =
\begin{cases}
\frac{1-\sqrt{1-4/P}}{4}    & P >\frac{25}{6}\\
\frac{1-3\sqrt{1-4/P}}{2}  & 4<P<\frac{25}{6}
\end{cases}
\end{equation}
The exponent is a decreasing function of $P$, see Fig.~\ref{fig:alpha}. When $P=\frac{25}{6}$, we have $m_3=C M^\frac{1}{5}$ with unknown amplitude $C$ and $m_2=\frac{5}{3}$, so that Eq.~\eqref{U-mass} becomes 
\begin{equation}
\frac{25}{3}M\,\frac{d W}{dM} = C M^\frac{1}{5} + \frac{5}{3} (W-1)
\end{equation}
giving $W\sim M^\frac{1}{5}\log M$. This solution together with \eqref{alpha:def} and \eqref{alpha:P} lead to \eqref{mass:uni} for $P>4$. 

In the percolating regime ($P<4$), the size distribution of unicycles has an exponential tail [see Eq.~\eqref{Us-star}]. Therefore, both $U$ and $W>U$ remain finite, although they are extremely large near the percolation transition. The average total number of unicycles is $U=C_-$. The exact formula \eqref{Cminus} yields $U \simeq \frac{25}{24} g^{-1}$ near the percolation transition, from which $U\sim e^{\pi/\sqrt{1-P/4}}$. The average mass of unicycles grows much faster when $P\uparrow 4$, apparently as $W\sim e^{2\pi/\sqrt{1-P/4}}$.

\end{document}